\documentclass[11pt,a4paper]{amsart}
\usepackage[utf8]{inputenc}
\usepackage[english]{babel}
\usepackage{amsmath}
\usepackage{amsfonts}
\usepackage{amssymb}
\usepackage{graphicx}
\usepackage{lmodern}
\usepackage[left=3.5cm,right=3.5cm,top=2.5cm,bottom=2.5cm]{geometry}

\makeindex             

\newcommand{\R}{\mathbb{R}} 
\newcommand{\N}{\mathbb{N}} 
\newcommand{\C}{\mathcal{C}} 
\newcommand{\PS}{\mathcal{P}} 
\DeclareMathOperator{\argmin}{argmin}
\newcommand{\e}{\epsilon} 
\newcommand{\norm}[1]{\left\lVert#1\right\rVert}
\DeclareMathOperator{\essInf}{ess inf}

\newcommand{\Hxp}{\text{H}(f(X^i(t)) - f(p^i))}
\newcommand{\Hvp}{\text{H}(f(v_f) - f(p^i))}
\newcommand{\Hxv}{\text{H}(f(X^i(t)) - f(v_f))}
\newcommand{\Hpv}{\text{H}(f(p^i)-f(v_f))}

\newtheorem{theorem}{Theorem}
\theoremstyle{plain}
\newtheorem{assumption}{Assumption}
\newtheorem{proposition}{Proposition}
\newtheorem{remark}{Remark}
\newtheorem{definition}{Definition}

\title{Trends in Consensus-based optimization}
\author[C. Totzeck]{Claudia Totzeck}
\address{University of Mannheim, B6, 68159 Mannheim }
\email{totzeck@uni-mannheim.de}

\begin{document}

\maketitle

\begin{abstract}In this chapter we give an overview of the consensus-based global optimization algorithm and its recent variants. We recall the formulation and analytical results of the original model, then we discuss variants using component-wise independent or common noise. In combination with mini-batch approaches those variants were tailored for  machine learning applications. Moreover, it turns out that the analytical estimates are dimension independent, which is useful for high-dimensional problems. We discuss the relationship of consensus-based optimization with particle swarm optimization, a method widely used in the engineering community. Then we survey a variant of consensus-based optimization that is proposed for global optimization problems constrained to hyper-surfaces. We conclude the chapter with remarks on applications, preprints and open problems.
\end{abstract}

\section{Introduction}
\label{sec:1}
\begin{figure}
	\includegraphics[scale=0.3]{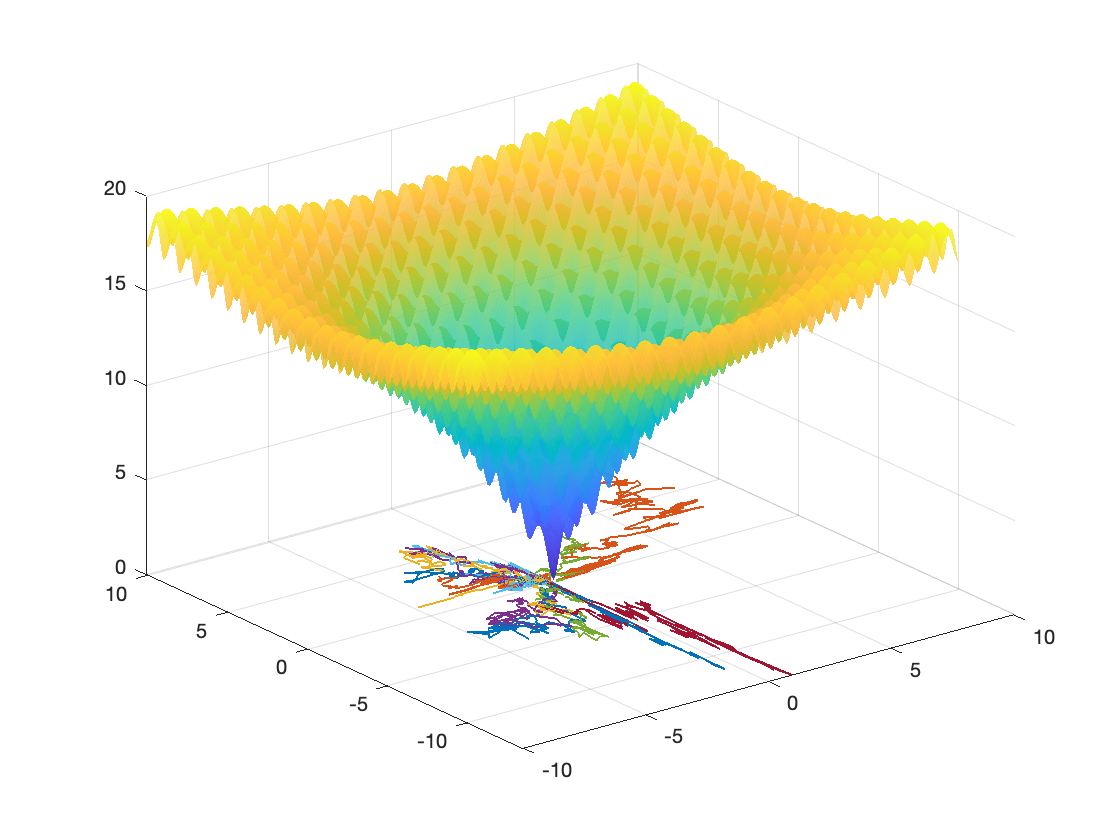}
	\caption{Plot of the Ackley \cite{benchmarks} benchmark function for global optimization in two dimensions with trajectories of one realisation of \eqref{eq:SDEcomponentwise} with $20$ particles visualized in the $xy$-plane.}
	\label{fig:ackley_trajectories}
\end{figure}

Global optimization tasks arise in various fields such as economics, finance, physics, clustering and artificial intelligence. In the most general form, these read
$$ \min\limits_{x\in \mathcal X} f(x) $$
for a given objective function $f$ and state space $\mathcal X$. Despite its simple description, the problem is nontrivial for nonconvex $f$ with possibly many local minima or constraint state spaces $X$, see Figure~\ref{fig:ackley_trajectories}. Its importance in various disciplines attracted the attention of many researchers to seek for solution strategies. Here, we focus on agent-based methods: on the one hand, there are biologically inspired methods as the ant colony optimization \cite{ant}, artificial bee colony optimization \cite{ABC} or Firefly Optimization \cite{Firefly}.

 On the other hand, wind driven optimization (WDO) \cite{WDO} is physically inspired as it models weather phenomena such as pressure and wind. The most popular agent-based global optimization algorithms are the Particle Swarm Optimization (PSO) and Simulated Annealing (SA). In PSO agents explore the state space while encountering a randomized drift towards the global best position seen by all the agents and a second drift towards their personal best positions. We will see more details on PSO below in Section~\ref{sec:PSO} where similarities and differences of CBO and PSO are discussed.
 SA is physically inspired, again, agents explore the state space. They are driven by noise terms that are diminishing as time evolves. The decrease of stochastic influence is called cool down and the particles are expected to concentrate at the best position seen by the particles during the exploration phase.

Most of the global optimization approaches are heuristics that have proven to give useful results in applications, but lack a rigorous analysis. Some proofs of convergence exist for SA in the context of image restoration and global optimization. These are mostly in the discrete setting and based on Markov Chains, see the survey \cite{convergenceSA} for more details.
 
A main objective in the modelling of the CBO scheme was to treat all particles identically, in particular, to circumvent the selection of a current best particle. In this way, one expects to have a corresponding mean-field scheme that can be utilized for the convergence analysis. Having this in mind, the CBO method is proposed as system of Stochastic Differential Equations (SDE) that mimics interacting agents communicating over a weighted mean. By construction the particles are expected to build a consensus at the position of the weighted mean that is located near the global minimizer of the functional.

To achieve this behaviour CBO combines ideas of swarm intelligence \cite{PSO} 
with approaches from consensus formation \cite{consensusFormation} 
in order to obtain a scheme that minimizes the objective function. CBO was first introduced in \cite{CBO1}, where formal relations to the mean-field equation and promising numerical results were shown. The main feature of the CBO algorithm is a weighted mean, $v_f.$  Particles with small function values have more influence in the weighted mean than particles with large function values. In this way the weighted mean is expected to be a good approximation of the global minimizer. All particles are driven by two terms. A drift term forcing them to move towards the weighted mean and a scaled diffusion allowing for exploration. In fact, whenever a particle is far away from the weighted mean, it explores its surroundings and tries to find a better position than the weighed average has. The scaling of the diffusion depends on the distance of the particle to the weighted mean. If the two coincide, the diffusion vanishes. Hence, the scheme allows for concentration at the position of the weighted mean.

The fact that the global minimizer is approximated with the help of the weighted mean is crucial when it comes to the mean-field limit. In particular, using the weighted mean the scheme circumvents to label any particle as current leader, or current global best, which would make the particles distinguishable and prevent us to carry out the mean-field limit. Formally, the limiting equation for `number of particles to infinity` is the PDE corresponding to the McKean process resulting from It\^{o} calculus applied to the SDE system \cite{CBO1,CBO2}. 
In \cite{CBO2} a rigorous analysis of the PDE method is performed. In particular, it is shown that the method converges to the minimizer of the global optimization scheme under some appropriate conditions. 

Another advantage of the communication via the weighed mean is a reduction of the computational effort. In fact, the communication with the weighted mean is of order $\mathcal O(N)$ for $N$ particles in the swarm. In other consensus algorithms each particle communicates with all other particle separately, leading to an effort of $\mathcal O(N^2)$ which suffers the curse of dimensionality when the swarm size grows. 

Recently, variants and extensions were proposed to improve the CBO method. Some approaches aim to enhance the performance in high-dimensional problems such as the ones arising in machine learning. Others extend the class of problems to be solved with CBO, for example, they allow for constrained state space $\mathcal X$.

In this survey we discuss these advances and compare them to the original method. The main part covers models which have been approved by peer-review. At the end we shed some light on recent preprints as well. Before we go into the details we shortly describe the ideas covered in the following. 

In \cite{CarrillJinZhu} the diffusion term was replaced by a component-wise diffusion, leading to a scheme that is robust with respect to the dimension of the state space. Indeed, the authors were able to show that many of the estimates shown in \cite{CBO2} hold without dimension-dependence for the scheme with component-wise diffusion. Moreover, the article introduces a mini-batch idea for the computation of the weighted mean. This reduces the computational cost and has positive effects on the performance in high dimensional scenarios. More details on this variant are discussed in Section~\ref{sec:Variant1}.
The authors of \cite{HaJinKimConvergence}
 replace the component-wise independent noise of the above variant by a component-wise common noise. This adaption facilitates the analysis of the scheme on the particle level. In fact, the authors show convergence of the variant directly on the particle level in contrast to \cite{CBO2,CarrillJinZhu},
  where the PDE formulation is employed for the analysis.
A variant that incorporates global in time information in order to approximate the personal best position seen by each of the particles is proposed in \cite{TotzeckWolfram}. It is shown that this variant is robust even if the initial distribution of particles is inconvenient. We discuss the scheme with global in time information and its relationship to PSO in Section~\ref{sec:PSO}.

In addition, there are variants that take care of optimization problems on constrained sets.  Box constraints are rather simple to handle. Dynamics constrained to hyper-surfaces, for example the sphere, need more sophisticated ideas \cite{Fornasier1}. We discuss  approaches for constraint sets in Section~\ref{sec:constraints} and in Section~\ref{sec:applications} 
we briefly comment on the performance of the CBO variants. For example in \cite{CarrillJinZhu} are comparisons to stochastic gradient descent (SGD) methods and several studies for global optimization benchmarks reported. We conclude with an outlook to future work and open problems.

\subsection{Notation and assumptions}
Let us first fix the notation and assumptions that are consistently used in the following sections. This has the advantage that the sections are self-consistent and one might jump to the variant of most interest right after the introduction.

We denote the dimension of the state space by $d \in \N$ and $N \in \N$ is the number of agents or particles in the swarm. The two notions, agents and particles, are used equally throughout the text. The state of the $i$-th agent is given by a vector $X^i \colon [0,T] \rightarrow \R^d, i=1,\dots,N,$ and we collect the states of all agents at time $t\in[0,T]$ in the vector $X(t) = (X^1(t), \dots, X^N(t)) \in \R^{dN}.$ The initial condition of the particles is denoted by $X_0^i \in\R^d$ for $i=1,\dots,N$ and we assume that $X_0^i$ are independent and identically distributed with $\text{law}(X^i_0) = \rho_0 \in \mathcal P(\R^d).$ The constants $\lambda, \sigma\ge0$ denote the drift and diffusion parameters, respectively. Some schemes incorporate a Heaviside function $H$ or a regularization $H^\epsilon$ thereof, which we fix as
\[
H(x) = \begin{cases} 1, & \text{for } x \ge 0, \\ 0, & \text{else} \end{cases}, \qquad \qquad H^\epsilon(x) = \frac12 + \frac12 \tanh\left( \frac x\epsilon \right).
\]
Moreover, we denote by $W^i,i=1,\dots,N$ independent $d$-dimensional Brownian motions. We consider the minimization problem
\begin{equation}\label{eq:GOP} \tag{P}
\min\limits_{x\in \mathcal X} f(x), 
\end{equation}
where  $f \colon \R^d \rightarrow \R_{\ge0}$ is a continuous function that admits a unique global minimizer $X_*\in\R^d$ and $\mathcal X = \R^d$ except for Section \ref{sec:variant5}, where we discuss state constraints and minimize $f$ on some hyper-surface $\Gamma \subset \R^d.$

\subsubsection{The weighted average}
As mentioned above, a weighted average or weighted mean plays a crucial role in all variants of CBO. For simplicity, we fix the weight function to be 
\begin{equation}\label{eq:weight}
	\omega_\alpha^f(x) =\exp(-\alpha f(x))
\end{equation}
throughout this review. Other choices are possible as well, but the weight function should be tailored to represent the task of finding a global minimum. 

Unless otherwise stated, the notion \textit{weighted mean} refers to the vector 
\begin{equation}\label{eq:vfParticle}
v_f = \frac{1 }{\sum_{i=1}^N \omega_\alpha^f(X^i(t))} \sum_{i=1}^N X^i(t) \omega_\alpha^f(X^i(t)) .
\end{equation}
Note that the objective function enters into the weight. Hence, due to \eqref{eq:weight} agents at locations with lower function values have more weight in the mean than agents located at positions with high function values. The parameter $\alpha$ controls this separation effect. Indeed, for $\alpha = 0$ all particles have the same weight and for $\alpha \rightarrow \infty$ we expect $v_f$ to approximate the global best of the agents, i.e.,
\[
v_f \approx \argmin_{i=1,\dots,N} f(X^i(t)).
\]
Note that the $ß\argmin$ may be set-valued in general. For simplicity, we assumed above that $f$ attains a unique minimizer. 

The argument for $\alpha \rightarrow \infty$ is related to the Laplace principle from large deviation theory \cite{DemboZeitouni}. In fact, under the assumption that the processes $X^i(t)$ are independent, we formally pass to the limit $N \rightarrow \infty$ to obtain
\[
\frac{1 }{\sum_{i=1}^N \omega_\alpha^f(X^i(t))} \sum_{i=1}^N X^i(t) \omega_\alpha^f(X^i(t))  \rightarrow \frac{1 }{\int \omega_\alpha^f(x) d\rho_t} \int x \omega_\alpha^f(x) d\rho_t 
\]
in distributional sense, with $\rho_t \in \mathcal P(\R^d)$ being the  Borel probability measure describing the one-particle mean-field distribution, which is assumed to be absolutely continuous w.r.t. the Lebesgue measure $dx$. Then, by Laplace principle \cite{CBO1} we have

\begin{proposition}\label{prop:LaplacePrinciple}
	Assume that $f \in \mathcal C_b(\R^d,\R), f \ge 0$, attains a unique global minimum at the point $X_* \in \R^d$ and let $\rho \in \mathcal P^{ac}(\R^d).$ Then, it holds
	\[
	\lim\limits_{\alpha \rightarrow \infty} \left(-\frac1\alpha \log \left( \int_{\R^d} e^{-\alpha f} d\rho \right)\right) = f(X_*).
	\]
\end{proposition}	
This property is the main motivation to choose the $\omega_\alpha^f$ as given in \eqref{eq:weight}. Note that uniqueness of the minimizer plays a role. If there were several global minimizers, the weighted mean would be in the convex hull of these and, in general, have a greater function value.

In the following section we recall the original statement of the CBO scheme, then we discuss recent variants. Readers familiar with the original scheme may jump directly to the variant of their interest.

\section{Consensus-based global optimization methods}
We begin this section with the original method as proposed in \cite{CBO1} and analysed in \cite{CBO2}. Then we move on to recent variants that were tailored to improve the method for high dimensional applications as arising in machine learning. The variants replace the diffusion term with either component-wise independent or component-wise common diffusion.

\subsection{Original statement of the method}
\label{sec:2}
The ideas behind and main features of CBO \cite{CBO1} are explained on the particle level. Then we formally pass to the mean-field level and review analytical results that discuss the formation of consensus near the global minimizer \cite{CBO2}.
\subsubsection{Particle scheme}
Consensus-based optimization was first introduced in \cite{CBO1} as a swarm dynamic that consists of $N$ coupled stochastic differential equations (SDE). The equation of the $i$-th agent is given by
\begin{equation}\label{eq:originalSDE}
dX^i(t) = -\lambda (X^i(t) - v_f) H^\epsilon(f(X^i(t))- f(v_f) ) dt + \sqrt{2} \sigma |X^i(t) - v_f| dW^i(t), 
\end{equation}
for $i=1,\dots,N$ and supplemented with initial data $X(0) = X_0.$ The system is coupled by the weighted average, $v_f,$ which appears in the equation of every agent. The first term on the right-hand side models a drift towards $v_f$. The greater the distance of the agent's position to $v_f$ the stronger the drift. The Heaviside function assures that the particle only moves towards $v_f$, if the function value of $v_f$ is better, i.e., smaller then the function value of the particle. The idea behind the diffusion term is similar. The diffusivity is scaled with the distance of $|X^i(t) - v_f|$, an agent far away from $v_f$ is allowed to explore its neighbourhood and possibly find a better position than $v_f$. While an agent close to $v_f$ is less diffusive and tends to keep its position. In particular, the  diffusion of particle $i$ vanishes if $X^i = v_f.$ This allows for concentration of the particles at $v_f$. 

\begin{remark}\label{rem:advantagesCBO} Let us emphasize some advantages of this dynamic:
	\begin{enumerate}
	\item \textbf{Indistinguishable particles:} Compared to other swarm intelligence schemes the dynamic does not depend on $\argmin_{X^i} f(X^i),$ but only of its approximation $v_f.$ Therefore, we may formally derive a limiting equation in mean-field sense as $N \rightarrow \infty,$ compare Section~\ref{sec:mean-field}, and use the PDE for the analytical investigation.
	\item \textbf{Interaction scales with $N$:} The coupling via $v_f$ has a huge advantage as well from a numerical point of view as we do not have binary interactions. The effort for the interaction of the agents scales only linearly in $N$. This is in contrast to many interaction models for crowd dynamics, where agents interact with all other agents at the same time, leading to a convolution term of order $\mathcal O(N^2).$
	\item \textbf{Exploration of full space:} Due to the term $|X^i(t) - v_f| dW^i(t)$, exploration takes place in $\R^d$ even if the $X^i$ are initial spanning only a subspace of $\R^d.$ This has a positive effect on the exploration if $N\ll d.$ 
 	\end{enumerate}
\end{remark}

\subsubsection*{Heaviside function}

In the original model, the Heaviside function was imposed to make concentration in local minima less probable. As reported in \cite{CBO1}, the deterministic scheme, $\sigma=0,$ with Heaviside function allows for a stationary solutions consisting of several Dirac measures located at level sets of $f$. For $\sigma > 0$ these solutions have probability zero, due to the Brownian motion. Moreover, it turned out in numerical studies that the scheme works fine without the multiplication of the Heaviside function and for the analytical investigation in \cite{CBO2} it was neglected. Therefore, we mainly focus on the scheme without Heaviside function in the following.

\subsubsection{Mean-field limit}\label{sec:mean-field}
Properties of the scheme were investigated on the mean-field level. Up to the author's knowledge there is no rigorous proof of the limiting equation so far. We therefore have to assume that propagation of chaos holds in order to derive the mean-field equation formally.

Let us assume that the \textit{propagation of chaos} property holds, that means the distribution of all agents $X, \nu_t^N,$ satisfies $\nu_t^N \approx \rho_t^{\otimes N}, N \gg 1$ and therefore $X^i(t)$ are approximately independently $\rho_t$-distributed. Then
\[
\frac1N \sum_{i=1}^N \omega_\alpha^f(X^i(t)) \approx \int_{R^d} \omega_\alpha^f(x) d\rho_t, \qquad \frac{1}{N} \sum_{i=1}^N X^i(t) w_\alpha^f(X^i(t)) \approx \int_{\R^d} x \omega_\alpha^f(x) d\rho_t, 
\]
due to the law of large numbers. Hence, $v_f \approx v_f[\rho_t]$ and we obtain the \textit{McKean nonlinear process} 
\begin{subequations}
	\begin{align}\label{eq:McKean}
		d\bar X(t) = -\lambda (\bar X(t)-v_f[\rho_t])\,d t + \sqrt{2} \sigma |\bar X(t)-v_f[\rho_t]| d W_t,
	\end{align}
	where the weighted average reads
	\begin{align}\label{eq:vfmeanfield}
		v_f[\rho_t ] = \frac{1}{\int_{\R^d} \omega_f^\alpha d\rho_t} \int_{\R^d} x\,\omega_f^\alpha d\rho_t,\qquad \rho_t = \text{law}(\bar X(t)).
	\end{align}
\end{subequations}
Equation \eqref{eq:McKean} may be equivalently expressed as the Fokker--Planck equation
\begin{subequations}\label{eq:model-pde}
\begin{gather}
	\partial_t \rho_t = \Delta (\kappa[\rho_t]\rho_t) + \text{div}(\mu[\rho_t]\rho_t), \\
	\kappa[\rho_t](x) = \sigma^2|x-v_f[\rho_t]|^2,\qquad \mu[\rho_t](x) = -\lambda(x-v_f[\rho_t]),
\end{gather}
\end{subequations}
which describes the evolution of the law corresponding to the Mc-Kean nonlinear process $\{\bar X(t)\in\R^d\,|\,t\ge 0\}$. 

The presence of $v_f$ makes the Fokker--Planck equation nonlinear and nonlocal in both the drift and diffusion part. This is nonstandard in the literature and raised several analytical and numerical questions that were addressed in \cite{CBO2}. We recall the main results in the following.

\subsubsection{Analytical results for the original scheme without Heaviside function}

The first statement considers the well-posedness of the particle dynamic, see Theorem~2.1 in \cite{CBO2} for the proof.
\begin{theorem}
Let the objective function $f$ be locally Lipschitz continuous.  For every $N \in \N$ system \eqref{eq:originalSDE} has a unique strong solution $\{ X_t^{N} \colon t \ge 0\}$  for any initial condition $X_0^{(N)}$ satisfying $\mathbb E |X_0^{(N)}|^2 < \infty. $
\end{theorem}

For the original particle scheme there is neither a proof for consensus formation nor for convergence to the global minimizer. These kinds of results were only addressed on the mean-field level after a formal limiting procedure as $N \rightarrow \infty.$ A rigorous proof of this limit is open up to the author's knowledge. The following results and some first estimates in the direction of a rigorous proof of the mean-field limit are reported in \cite{CBO2}. 

The well-posedness of the mean-field equation is established for two classes of objective functions. One result considers only bounded objective functions and another result is for objective functions with quadratic growth at infinity. Both versions are based on the following assumption:

\begin{assumption}\label{ass} To obtain the well-posedness results of the mean-field equation we assume that it holds:
	\begin{enumerate}
		\item The cost function $f\colon\R^d\to\R$ is bounded from below with $\underline f := \inf f$.
		\item There exist constants $L_f$ and $c_u>0$ such that
		\begin{align}
			\left\lbrace\quad\begin{aligned}
				|f(x)-f(y)| &\le L_f(|x|+|y|)|x-y| \quad \text{for all\, $x,y\in\R^d$},\\
				f(x)-\underline f &\le c_u(1+|x|^2)\quad \text{for all\, $x\in\R^d$}.
			\end{aligned}\right.\tag{A1}\label{ass:1}
		\end{align}
	\end{enumerate}
\end{assumption}

\begin{definition}\label{def:quadraticGrowth}
We say that a function has \textit{quadratic growth} if there exist constants $M>0$ and $c_l>0$ such that
\begin{align}
	f(x)- \underline f \ge c_l|x|^2\quad\text{for\, $|x|\ge M$}\tag{A2}\label{ass:2}.
\end{align}
\end{definition}

\begin{theorem}\label{thm:well_posed_meanfield}
	Let $f$ be bounded or have quadratic growth, let Assumption~\ref{ass} hold and $\rho_0\in\PS_4(\R^d)$.
	Then there exists a unique nonlinear process $\bar X\in \C([0,T],\R^d)$, $T>0$, satisfying 
	\begin{align*}
		d\bar X_t = -\lambda(\bar X_t - v_f[\rho_t])\,dt + \sigma|\bar X_t - v_f[\rho_t]|dW_t,\qquad \rho_t=\text{law}(\bar X_t),
	\end{align*}
	in the strong sense, and $\rho\in \C([0,T],\PS_2(\R^d))$ satisfies the corresponding Fokker--Planck equation \eqref{eq:model-pde} in the weak sense with $\lim\nolimits_{t\to 0}\rho_t=\rho_0\in\PS_2(\R^d)$.
\end{theorem}
Both proofs are based on Schauders fixed-point argument and can be found in \cite{CBO2}. The main difference of the two versions is the argument for the bound of the second moment. This bound is needed in order to apply Gronwall theorem and to close the Schauder argument. 

Convergence of the scheme towards the global minimizer of the objective function is shown in two steps. The first step assures only the consensus formation. The second one shows that for appropriate parameter choices the consensus location is positioned near the global minimizer. Both results are of asymptotic nature. The consensus formation occurs for $t\rightarrow\infty$ and the approximation of the global minimizer depends on the choice of the weight parameter $\alpha$. For $\alpha\rightarrow\infty$ the location of consensus tends towards the global minimizer. 
For the concentration result we need this assumption.
\begin{assumption}\label{ass:3}
	We assume that $f \in \mathcal C^2(\R^d)$ satisfies additionally
	\begin{enumerate}
		\item $\inf f > 0$.
		\item $\|\nabla^2 f\|_\infty\le c_f$ and there exists constants $c_0,c_1>0$, such that
		\[
		\Delta f \le c_0 + c_1|\nabla f|^2\qquad\text{in\; $\R^d$}.
		\]
	\end{enumerate}
\end{assumption}
To show the concentration, we investigate the expectation and the variance of the density which are defined by
\[
E(\rho_t) = \int_{\mathcal X} x d\rho_t \quad \text{and}\quad V(\rho_t) = \frac12\int_{\mathcal X} |x - E(\rho_t)|^2 d\rho_t.
\]
 The details of the concentration procedure are given in \cite{CBO2}.
\begin{theorem}\label{thm:concentration}
	Let $f$ satisfy Assumption~\ref{ass:3} and let the parameters $\alpha,\lambda$ and $\sigma$ satisfy
	\[
	2\alpha e^{-2\alpha\underline{f}} (c_0 \sigma^2 + 2\lambda c_f) < \frac{3}{4},\qquad 2\lambda b_0^2 - K - 2d\sigma^2 b_0e^{-\alpha \underline f} \ge 0, 
	\]
	with $K=V(\rho_0)$ and $b_0=\|\omega^\alpha_f\|_{L^1(\rho_0)}$. Then $V(\rho_t)\le V(\rho_0)e^{-q t}$ with
	\[
	q = 2\big(\lambda - (d\sigma^2/b_0)e^{-\alpha \underline f}\big) \ge K/b_0^2.
	\]
	In particular, there exists a point $\tilde x\in\R^d$ for which $E(\rho_t)\to \tilde x$ and $v_f[\rho_t]\to \tilde x$ as $t\to\infty$.
\end{theorem}
So far, we just know that the density will concentrate at some point, $\tilde x$, the location of this point remains unknown. Finally, the following result assures that the concentration takes place in a neighbourhood of the global minimizer for appropriately chosen parameters.
\begin{theorem}\label{thm:approximate}
	Let $f$ satisfy Assumption~\ref{ass:3}. For any given $0<\e_0\ll1$ arbitrarily small, there exist some $\alpha_0\gg 1$ and appropriate parameters $(\lambda,\sigma)$ such that uniform consensus is obtained at a point $\tilde x\in B_{\e_0}(x_*)$. More precisely, we have that $\rho_t \to \delta_{\tilde x}$ for $t\to\infty$, with $\tilde x\in B_{\e_0}(x_*)$.
\end{theorem}

With this result we conclude the survey of the analytical results on the original scheme of consensus-based global optimization proposed in \cite{CBO1} and analysed in \cite{CBO2}. 

\subsubsection{Numerical methods}
It is important to notice that the success of the CBO method is far more dependent on the function evaluations than on the accuracy of the numerical scheme. In fact, whenever a particle hits the global minimum of the function, the weighted average $v_f$ is assumed to move to this position and then concentration takes place. 

Having this in mind, most of the numerical simulations use basic algorithms such as the Euler-Maruyama scheme \cite{EulerMaruyama}.

In \cite{CBO1} the formal mean-field limit is underlined by the comparison of numerical results on the particle level with the solution of the candidate equation on the mean-field level. The PDE is solved with the help of a discontinuous Galerkin approach in combination with a Strang splitting. The convective part is solved with the local Lax-Friedrichs scheme and the diffusion part semi-implicitly.

In the following we discuss variants of this method that aim to enhance the performance or extend the class of optimization problems admissible for CBO. We begin with a variant that appears like a slight modification of the above algorithm. However, it has a major impact on the convergence results, especially in high dimensions.

\subsection{Variant 1: Component-wise diffusion and random batches}\label{sec:Variant1}
At first glance, the variant with component-wise independent noise in \cite{CarrillJinZhu} seems to be a minor modification of the original dynamic. Nevertheless, it turns out that the estimates of the convergence results become independent of the dimension of the state space. This is an advantage, especially when the method is considered for high-dimensional problems, for example, arising in machine learning. In addition to the component-wise diffusion, the authors propose to use mini-batches, a popular approach in many stochastic gradient descent methods \cite{MiniBatch}.

\subsubsection{Component-wise geometric Brownian motion}
The dynamic with component-wise geometric Brownian motion reads
\begin{equation}\label{eq:SDEcomponentwise}
	dX^i(t) = -\lambda (X^i(t) - v_f) dt + \hat \sigma \sum_{k=1}^d (X^i(t) - v_f)_k dW_k^i(t) \vec{e}_k, 
\end{equation}
for $i=1,\dots,N$ and is supplemented with initial data $X(0) = X_0.$ Here $\vec e_k$ denotes the $k$-th unit vector in $\R^d$, $(X^i(t) - v_f)_k$ is the $k$-th entry of the difference and  $W_k^i$ are independent standard Brownian motions. The weighted mean, $v_f,$ is given in \eqref{eq:vfParticle}.

\begin{remark}\label{rem:CBOcomp} Let us mention two differences between \eqref{eq:originalSDE} and \eqref{eq:SDEcomponentwise}:
	\begin{enumerate}
	\item \textbf{Component-wise noise:} The component-wise diffusion in \eqref{eq:SDEcomponentwise} scales the distance of $X^i_k$ and $v_f$ element-wise. In case one component of the two coincide, this component of $X^i$ does not change.
	\item \textbf{Diffusion constants:} The slight difference between the diffusion constants in \eqref{eq:originalSDE} and \eqref{eq:SDEcomponentwise} $\hat \sigma = \sqrt{2} \sigma$ has no significant influence on the performance of the scheme. 
	\end{enumerate}
\end{remark}

The aforementioned dimension-independence of the component-wise diffusion can be seen with the help of a simple computation \cite{CarrillJinZhu}. Let us fix the weighed average $v_f$ at an arbitrary position $a$. Then, for the dynamics in \eqref{eq:originalSDE} we find
\[
\frac{d}{dt} \mathbb E |X(t) - a|^2 = -2\lambda \mathbb E |X(t) - a|^2 + \sigma^2 \sum_{i=1}^d \mathbb E |X(t) - a|^2 = (-2\lambda + \sigma^2 d) \mathbb E |X(t) - a|^2.
\]
This investigation of the second moment shows, that concentration occurs whenever the condition $2\lambda > \sigma^2 d$ is satisfied.
In contrast, the same computation for \eqref{eq:SDEcomponentwise} yields
\[
\frac{d}{dt} \mathbb E |X(t) - a|^2 = -2\lambda \mathbb E |X(t) - a|^2 + \sigma^2 \sum_{i=1}^d \mathbb E |X(t) - a|^2_i = (-2\lambda + \sigma^2 ) \mathbb E |X(t) - a|^2.
\]
The condition for concentration changes to $2\lambda > \sigma^2$. In particular, it is independent of the dimension $d.$ 

It can be proven that all estimates needed for the analysis of well-posedness, concentration and convergence towards the global minimizer on the mean-field level are independent of the dimension for the component-wise diffusion variant. Instead of rewriting the statements here, we refer to \cite{CarrillJinZhu} for all details and proceed with the second interesting feature proposed in the article.
 
 \subsubsection{Random batch method}  
 The second novelty proposed in \cite{CarrillJinZhu} is to apply the random-mini batch strategy \cite{JinLiLiu} in two levels: first, instead of evaluating $f(X^i(t))$ for every particle $i=1,\dots,N$ in every time step, $q$  random subsets $J^\theta\subset \{1,\dots,N\}$ with size $|J^\theta|=M \ll N$ and $\theta = 1,\dots,q$ are drawn and for each of them an empirical expectation $\hat f(X^\theta)$ is computed. Based on these function evaluations a weighted mean is calculated for every batch. Now, one can choose to update the positions of particles by \eqref{eq:SDEcomponentwise} only for the particle in the batch, or apply the dynamics to all $N$ particles. 
 For simplicity, we present a version of the algorithm in \cite{CarrillJinZhu} adapted to the general problem \eqref{eq:GOP}.  Note that there is an additional parameter, $\gamma_{k,\theta},$ called \textit{learning rate} following the machine learning terminology. \\

\noindent\textbf{Algorithm 1}\newline
	Generate $\{X_0^i \in \R^d\}_{i=1}^N$ according to the same distribution $\rho_0$.   
	Set the remainder set $\mathcal{R}_0$ to be empty. For $k=0,1,2, \dots $, do the following:
	\begin{itemize}
	\item [-] Concatenate $\mathcal{R}_k$ and a random permutation $\mathcal{P}_k$ of the indices $\{1, 2,\cdots, N\}$ to form a list $\mathcal{I}_k=[\mathcal{R}_k, \mathcal{P}_k]$.  Pick $q=\lfloor \frac{N+|\mathcal{R}_k|}{M}\rfloor$ sets of size $M\ll N$ from the list $\mathcal{I}_k$ in order to get batches $B_1^k, B_2^k, \dots, B_q^k$ and set the remaining indices to be $\mathcal{R}_{k+1}$. Here, $|\mathcal{R}_k|$ means the number of elements in $\mathcal{R}_k$.

	\item [-]  For each $B_{\theta}^k$ ($\theta=1,\cdots, q$), do the following
	
	\begin{enumerate}
		\item Calculate the function values (or approximated values) of $f$ at the location of the particles in $B_\theta^k$ by $f^j:= f(X^j),~\forall j\in B_\theta^k$.
		
		\item  Update $v_{k,\theta}$ according to the following weighted average, 
		\begin{equation*}
			v_{k,\theta} = \frac{1}{\sum_{j\in B_{\theta}^k} \mu_j}\sum_{j\in B_{\theta}^k} X^i \mu_j ,\quad \text{with}\quad  \mu_j = e^{-\alpha f^j}.
		\end{equation*}
		\item Update $X^j$ for $j\in \mathcal{J}_{k,\theta}$ as follows, 
		\begin{equation*}
			X^j \leftarrow X^j - \lambda \gamma_{k,\theta}( X^j - v_{k,\theta})+ \sigma_{k,\theta} \sqrt{\gamma_{k,\theta}} \sum_{i=1}^d \vec{e}_i(X^j - v_{k,\theta})_i,    z_i^j,~~z_i^j\sim\mathcal{N}(0, 1),
		\end{equation*}
		where $\gamma_{k,\theta}$ is chosen suitably and there are two options for $\mathcal{J}_{k,\theta}$:
		\begin{equation*}
			\begin{aligned}
				&\text{{\it partial updates: }}\quad \mathcal{J}_{k,\theta}=B_{\theta}^k,\quad \text{or} \quad
				&\text{{\it full updates: }}\quad  \mathcal{J}_{k,\theta}=\{1,\cdots, N\}.
			\end{aligned}    
		\end{equation*}
	\end{enumerate}	
	\item [-] Check the {\bf Stopping criterion:} 
	\begin{equation*}
		\frac{1}{d} \norm{ \Delta x}_2^2 \leq \epsilon,
	\end{equation*}
	where $\norm{ \cdot} _2$ is the Euclidean norm and $\Delta v$ is the difference between two most recent $v_{k,\theta}$.
	If this is not satisfied, repeat. 
	\end{itemize}

Note that due to the mini-batch evaluation additional noise is added to the algorithm. The authors discuss in \cite{CarrillJinZhu} that this additional noise causes the algorithm to work fine even without the geometric Brownian motion. For details and additional ideas on how to improve the convergence for objective functions with a typical machine learning structure we refer to \cite{CarrillJinZhu}.

We conclude this section with some ideas on the numerical implementation and the performance of the algorithm with random batches and component-wise geometric Brownian motion. 

\subsubsection{Implementation and numerical results}

A typical challenge is to avoid overshooting which refers to oscillations around $v$ in our context. The authors propose two approaches to do so.

 First, the drift and diffusion part of the scheme can be split. Then the drift part can be computed explicitly using 
\[
\hat X_k^j = v_{k} + (X_k^j - v_k) e^{-\lambda \gamma}
\]
which corresponds to a scheme for solving the ODE $dX^j = -\lambda(X^j - v)$ on the interval $t \in [k\gamma, k(\gamma+1)].$ The diffusion update is given by
\[
X_{k+1}^j = \hat X_k^j + \sigma\sqrt{\gamma} \sum_{i=1}^d \vec{e}_i \left( \hat X_k^j - v \right)_i z_i^j.
\]

Second, they propose to freeze the weighted average over fixed time intervals. On each of these intervals the geometric Brownian motion can be solved by
\[
X_{k+1}^j = v + \sum_{i=1}^d \vec{e}_i \left( \hat X_k^j - v \right)_i\, \exp((-\lambda - \frac12 \sigma^2) \gamma + \sigma \sqrt{\gamma}z_i^j).
\]
Moreover, they report that the splitting and the freezing approach lead to comparable results in most numerical simulations. For more details, see \cite{CarrillJinZhu}.

The aforementioned paper reports results of three numerical studies. The first is a proof of concept using a one-dimensional objective function with many local minima and oscillatory behaviour.  The second study compares the method to the performance a stochastic gradient descent method applied to the MNIST data set. Finally, results for a test function in high dimensions with many local minima are provided. 

The test cases show that the proposed CBO algorithm with component-wise Brownian motion and mini-batches outperforms the stochastic gradient descent algorithm. Moreover, it turns out that the version with mini-batches leads to better results than the one with full evaluations in case of the MNIST data set. For more detailed discussions and studies of the influence of $\alpha$ and $N$ on the performance we refer the reader to the original article \cite{CarrillJinZhu}. 

\subsection{Variant 2: Component-wise common diffusion}
The idea of component-wise diffusion plays a role as well in \cite{HaJinKimConvergence,HaJinKimErrorEstimates} with the main difference that the component-wise noise is \textit{common} for all particles, that means, the dynamic is given by
\begin{equation}\label{eq:SDEcomponentwiseCommon}
	dX^i(t) = -\lambda (X^i(t) - v_f) dt + \hat \sigma \sum_{k=1}^d (X^i(t) - v_f)_k dW_k(t) \vec{e}_k, 
\end{equation}
where $W_k$ are i.i.d.~one-dimensional Brownian motions. The dynamic is supplemented with initial conditions $X^i(0) = X_0^i$ and $v_f$ as above. Note that the Brownian motion does not depend on the specific particle $i$ therefore all particles encounter a common noise.

In addition to the continuous-time particle scheme given above, the article discusses a time-discrete version. Let $h >0$ denote the time step, i.e., $t = nh$ we set
$ X_n^i := X^i(nh).$ The discrete algorithm reads
\begin{align}\label{eq:commonNoiseDiscrete}
X_{n+1}^i = X_n^i -\lambda h(X_n^i - v_f) + \sigma \sqrt{h}  \sum_{k=1}^d (X_{n}^i - v_f)_k Z_n^k \vec{e}_k,
\end{align}
where $\{Z_n^k\}_{n,k}$ are i.i.d.~standard normal distributed random variables, $Z_n^k \sim \mathcal N(0,1).$ Note that compared to \cite{HaJinKimConvergence} the notation was adjusted for the sake of a consistent presentation.

\subsubsection{Analytical results}

The common noise approach has the advantage that a convergence study can be done directly on the level of particles without passing to the mean-field level. Similar to the strategy of the proof on the mean-field level the convergence proof for the common noise scheme is split into two parts: first, under certain conditions on the drift and diffusion parameters, a general convergence to consensus result for $t\to \infty$ is shown. In a second step the authors provide sufficient conditions on the system parameters and initial data which guarantee that the location of the consensus is in a small neighbourhood of the global minimum almost surely. The conditions on the parameters are independent of the dimension similar to Variant 1 (see Section \ref{sec:Variant1}).

Despite these two main results some properties of the continuous and discrete deterministic schemes are discussed. In fact, it is proven that the convex hull of the particles following the deterministic (both time-continuous and time-discrete) schemes are contractive as time evolves. The convergence to a  consensus state is a direct consequence. 

The same contraction property is not given for the scheme with noise. Nevertheless, for the common noise approach the relative difference of two particles satisfies a geometric Brownian motion. Hence, an exact solution can be established using stochastic calculus. This implies that the relative state difference converges almost surely. The details of the theorem are as follows.
\begin{theorem}
Let $X^i(t)$ be the $i$-the agent of a solution to \eqref{eq:SDEcomponentwiseCommon}. Then for $i\ne j = 1,\dots,N$  it holds
\[
\mathbb E |X^i(t) - X^j(t)|^2 = e^{-(2\lambda - \sigma^2)t} \mathbb E |X_0^i - X_0^j|^2, \quad t >0.
\]
In particular $L^2$-consensus emerges if and only if $\lambda - \frac{\sigma^2}{2} >0.$
\end{theorem}

A similar results is obtained for the time-discrete dynamic \eqref{eq:commonNoiseDiscrete}. The condition for the convergence depends on $\sigma, \lambda$ and $h.$ In fact, several different conditions are discussed. For details we refer to \cite{HaJinKimConvergence}.

The second step which shows that for well-chosen parameters the consensus state is located in a neighbourhood of the global minimizer is more involved. Here, we only state the main result which needs the following assumption.
\begin{assumption}\label{ass:common} We assume $f$ and the initial conditions satisfy:
\begin{enumerate}
\item $f \in C_b^2(\R^d)$ with $\inf\limits_{x\in\R^d} f(x) > 0$ and $$C_L := \max \Big\{ \sup\limits_{x\in\R^d} \| \nabla^2 f(x) \|_2 , \max\limits_{1\le l \le d} \sup\limits_{x\in\R^d} |\partial_l^2 f(x)|  \Big\} < \infty.$$
\item For some $\epsilon\in(0,1)$ the initial conditions $X_0^i$ are i.i.d.~with $X_0^i \sim X_\text{in}$ for some random variable $X_\text{in}$ which satisfies 
$$
(1-\epsilon)\mathbb E[e^{-\alpha f(X_\text{in})}] \ge \frac{2\lambda + \sigma^2}{2\lambda - \sigma^2} C_L \alpha e^{-\alpha f(X_*)} \sum_{l=1}^d \mathbb E \left[ \max\limits_{1 \le i \le N} (X_0^{i} - v_f(0))_l \right].
$$
\end{enumerate} 
\end{assumption}
\begin{theorem}Let Assumption \ref{ass:common} hold and suppose $2\lambda > \sigma^2.$ Then for a solution $X$ to \eqref{eq:SDEcomponentwiseCommon} it holds
$$ \lim_{t\rightarrow \infty}\essInf\limits_{\omega} f(X_t^i(\omega)) \le \essInf\limits_\omega f(X_\text{in}(\omega)) + E(\alpha) $$
for some function $E(\alpha)$ with $\lim\limits_{\alpha\to \infty} E(\alpha) = 0.$ In particular, if the global minimizer $X_*$ is contained in the support of $\text{law}(X_\text{in})=\rho_0$ then
$$  \lim_{t\rightarrow \infty}\essInf\limits_{\omega} f(X_t^i(\omega)) \le  f(X_*) + E(\alpha).  $$
\end{theorem}
The convergence of the time-discrete algorithm was not established in  \cite{HaJinKimConvergence}  due to the lack of a discrete analogue of It\^o's stochastic calculus. In a subsequent article \cite{HaJinKimErrorEstimates} the authors give an elementary convergence and error analysis for the time-discrete version \eqref{eq:commonNoiseDiscrete} under some additional regularity conditions on $f$. Moreover, exponential decay rates of the distances between the particles are established. The proofs are rather technical and go beyond the scope of this survey. We therefore refer the interested reader to \cite{HaJinKimErrorEstimates}.

\subsubsection{Numerical results}
A priori it is not clear how the common noise algorithm performs compared to the well-tested component-wise noise version in Section~\ref{sec:Variant1}. In \cite{HaJinKimConvergence} some numerical results of the common noise algorithm are provided. They underline the analytical results on the convergence of the distance of two particles and indicate that also the common noise version leads to reasonable results. A large-scale comparison of Variant 1 and the common noise scheme of this section is missing up to the author's knowledge.

\section{Relationship of CBO and Particle Swarm Optimization}\label{sec:PSO}
Consensus-based optimization is inspired by Particle Swarm Optimization (PSO) schemes \cite{PSO}. It is worthwhile to compare the methods to gain further insight to their behaviour, performance and the qualities. Let us recall the formulation of the PSO dynamic \cite{PoliKennedyBlackwell}: the update for the $i$-th particle is given by
\begin{align*}
X^i &\leftarrow X^i + V^i, \quad i=1,\dots,N \\
V^i &\leftarrow \omega V^i + \sum_{k=1}^d \Big(U_{1,k}^i (p_\text{personal} - X^i)_k + U_{2,k}^i (p_\text{global} - X^i)_k\Big), 
\end{align*}
where $U_{1/2}^i$ are $d$-dimensional vectors of random numbers which are uniformly distributed  in $[0,\phi_1]$ and $[0,\phi_2],$ respectively. $p_\text{global}$ denotes the best position that any of the particles has seen and $p_\text{personal}$ denotes the best position particle $i$ has seen. The parameters $\phi_{1,2}$ define the magnitude of the stochastic influences and $\omega$ can be interpreted as friction parameter. $V^i$ is originally kept within box constraints, given by the range $[-V_\text{max},V_\text{max}].$ 
In contrast to the first order dynamic of CBO, PSO is of second order which may lead to inertia effects.
 Moreover, the stochastic influence does not vanish, therefore one cannot expect any kind of consensus formation. The approximation of the global best is $p_\text{global}$ whenever the PSO algorithm is stopped. The global best information in PSO prevents a direct passage to the mean-field limit.

The main ingredient of CBO is the weighted average $v_f.$ For $\alpha \gg 1$ it can be interpreted as an approximation of the current best particle position.  Here, we use best in the sense that the function value is the lowest compared to the function values of all other particles. This current best particle does move only slightly, as its distance to $v_f$ is small and therefore the drift and diffusion terms are small. Therefore the current best particle can as well be interpreted as the global best position seen so far. Hence, $v_f$ is the analogue of $p_\text{global}$ in PSO. In addition, the PSO dynamic includes the so called local best position, which refers to the best position that each of the particles has seen. This local best is modelled in \cite{TotzeckWolfram} using a memory effect.   The same local best is mentioned as well in a recent preprint \cite{GrassiPareschi} which additionally considers a continuous description of PSO and computes the corresponding macroscopic equations to clarify the relationship of PSO and CBO. The details of the CBO with local or personal best information are given in the following.

\subsection{Variant 4: Personal best information}
Consensus-based optimization with global and local best in the sense of PSO is proposed in \cite{TotzeckWolfram} and based on the component-wise diffusion variant (see Section~\ref{sec:Variant1}). The dynamic reads as follows
	\begin{align}\label{eq:particlePSO}
	dX^i(t) &= \left[ -\lambda(t,X) (X^i(t)-v_f)- \mu(t,X) (X^i(t)-p^i(t) \right]\,d t \nonumber \\
	&\qquad\qquad+ \sqrt{2} \sigma \sum_{k=1}^d (X^i(t) - v_f)_k dW_k^i(t) \vec{e}_k, \qquad i=1,\dots,N
\end{align}
with $v_f$ as given above and the personal best is modelled by
$$ p^i(t) = \begin{cases} X_0^i, & t = 0, \\\int_0^t X^i(s) \exp(-\beta f(X^i(s))) ds \Big/ \int_0^t \exp(-\beta f(X^i(s))) ds, &\text{otherwise.} \end{cases}  $$
This personal best approximation uses the same idea as $v_f$ but with respect to time in contrast to the integral over the state space. Again by Laplace principle (see Proposition~\ref{prop:LaplacePrinciple}), we expect that $p^i(t)$ approximates the best position particle $X^i$ has seen up to time $t.$

\begin{remark}
To circumvent the integral over time, it is tempting to rewrite the numerator and denominator of $p^i$ as SDE. Notice that the initial condition of each personal best would need to be positioned at zero in order to obtain the exact definition above.
\end{remark}
To make sure that particles do not get stuck in the middle, each particle has to choose whether it moves towards $v_f$ or towards its personal best $p^i.$ As we aim for a global minimizer, we assume that this decision is based on the cost functional values, which motivates to set the prefactors $\lambda$ and $\mu$ as
\begin{align*} 
	\lambda(t,X) &=  \Hxv\, \Hpv,\\ \mu(t,X) &= \Hxp\, \Hvp.
\end{align*}
This is leading to the following behaviour:
\begin{enumerate}
	\item [-] If $f(v_f)$ is smaller than $f(X^i)$ and $f(p^i)$, the particle moves towards $v_f$.
	\item [-] If $f(p^i)$ is smaller than $f(X^i)$  and $f(v_f)$, the particle moves towards $p^i$.
	\item [-] If none of the above holds, the particle still explores the function landscape via Brownian motion until it reaches the global best $v_f$.
\end{enumerate}
Using a regularized version of the Heaviside function, $H^\epsilon,$ the well-posedness of the above system is proven in \cite{TotzeckWolfram}. There is no mean-field result and no convergence result reported. 

\subsubsection{Performance}
Note that the additional evaluation of the personal best position has minor impact on the computational costs as the time integrals in $p^i$ allow for an accumulative computation.
Note further, that even though the Heaviside function needs to be regularized for the analysis, the numerical results can work with the original Heaviside formulation.

The numerical results indicate that the personal best information raises the probability of finding the global best position, if few particles are involved in the search. As the number of particles needed for satisfying results, depends on the dimension of the state space, this result is particularly important in high dimensions. If the number of particles is large enough, no significant influence of the personal best information is noted.

\section{CBO with state constraints}\label{sec:constraints}
Many global optimization tasks have a constrained state space. The simplest version of constraints are box constraints. These can be included into each of the aforementioned CBO versions by projecting particles back into the box, whenever they are about to leave it. 

The situation is more complicated, when the state space is given in form of a hyper-surface of $\R^d$. For example the sphere 
\[
S^2 = \{ x \in \R^3 \colon |x| = 1  \}
\] is a hyper-surface of $\R^3.$ In \cite{Fornasier1, Fornasier2} a variant of CBO on such hyper-surfaces is proposed. The first paper is concerned with the well-posedness and the mean-field limit of the variant and the second article discusses the convergence to global minimizers and applications in machine learning. A major advantage of this variant is the fact, that compactness is assured by the constraint. Therefore the mean-field limit can be established rigorously. In the following we discuss the details \cite{Fornasier1}.

\subsection{Variant 5: Dynamics constrained to hyper-surfaces}\label{sec:variant5}
The restriction to the hyper-surface leads to a new formulation of the optimization problem
\[
\min\limits_{x \in \Gamma} f(x),
\]
where $\Gamma$ represents the hyper-surface and $f \colon \R^d \rightarrow \R$ as above. We assume that $\Gamma$ is a connected and smooth compact hyper-surface embedded in $\R^d$ which is represented as zero-level set of a signed distance function $\gamma$ with $|\gamma(x)| = \text{dist}(x,\Gamma)$ leading to 
\[
\Gamma = \{ x\in \R^d \colon \gamma(x) = 0 \}.
\] 
If $\partial \Gamma = \emptyset$ we assume for simplicity that $\gamma < 0$ on the interior of $\Gamma$ and $\gamma >0$ on the exterior. The gradient, $\nabla \gamma,$ is the outward unit normal on $\Gamma$ whenever $\gamma$ is defined. In addition, we assume that there exists an open neighbourhood $\hat \Gamma$ of $\Gamma$ such that $\gamma \in \C^3(\hat \Gamma).$ All these assumptions allow us to work with the Laplace-Beltrami operator. For example, for the sphere $\mathbb S^{d-1}$ we can choose $$\gamma(x) = |x|-1 \quad\text{with} \quad\nabla \gamma(x) = \frac{x}{|x|} \quad \text{and} \quad \Delta\gamma(x) = \frac{d-1}{|x|}.$$

In \cite{Fornasier1} a Kuramoto-Vicsek type dynamic is proposed as 
\begin{align}\label{CBOconstraint}
dX^i(t) &= -\lambda P(X^i(t))(X^i(t) - v_f) dt + \sigma |X^i(t) - v_f|P(X^i(t)) dB^i(t) \nonumber\\ &\quad\qquad-\frac{\sigma^2}{2}(X^i(t) - v_f)^2 \Delta \gamma(X^i(t)) \nabla \gamma(X^i(t)) dt, \quad i=1,\dots,N,
\end{align}
 with initial condition $X(0) = X_0.$ In contrast to the aforementioned schemes there appears a projection operator $P$ defined by
 $$ P(x) = I - \nabla \gamma(x) \nabla\gamma(x)^T. $$
 For the sphere we obtain the $P(x) = I - \frac{x\, x^T}{|x|^2}.$ In addition to this projection there appears a third term in  \eqref{CBOconstraint}. The two mechanisms ensure that the dynamics stays on the hyper-surface $\Gamma$.  
 
\begin{remark}
Note that the dynamic is described in $\R^d.$ On the one hand this allows for a simple statement of the scheme. On the other hand it is likely that the weighted average is not positioned at $\Gamma,$ i.e., $v_f \notin \Gamma.$ This is caused by the averaging of particles on a hyper-surface. Nevertheless, for $\alpha \gg 1,$ $v_f$ approximates the current best particle which is contained in $\Gamma$ due to the projection and correction terms.  
\end{remark}
 
 The constraint enables us to give rigorous arguments for the limit $N\to\infty,$ which results in the nonlocal, nonlinear Fokker-Planck equation
 \begin{align*}
\partial_t \rho_t = \lambda \nabla_\Gamma \cdot [P(x) (x - v_f) \rho_t] + \frac{\sigma^2}{2}\Delta_\Gamma (|x - v_f|^2 \rho_t), \quad t>0, x\in\Gamma,
\end{align*} 
with initial condition $\rho_0 \in \mathcal P(\Gamma).$ The operators $\nabla_\Gamma$ and $\Delta_\Gamma$ denote the divergence and Laplace-Beltrami operator on the hyper-surface $\Gamma,$ respectively.
In the following we summarize the analytical results for this variant which are reported in \cite{Fornasier1}.

\subsubsection{Analytical results}
The following analytical results focus on the well-posedness and the rigorous mean-field limit of the constrained scheme.

As the dynamic is living in $\R^d$ there are some technical issues with $P, \Delta\gamma$ and $\nabla \gamma.$ In fact, these are not defined for $x=0$ and the authors propose to replace them with regularizations. Moreover, a regularized extension of $f,$ called $\tilde f,$ is introduced.

\begin{assumption}\label{ass:constraint} Let $\tilde f$ be globally Lipschitz continuous and such that it holds: 
	\begin{gather*}
	\tilde f(x) = f(x) \text{ for } x \in \hat \Gamma, \\
	\tilde f(x) - \tilde f(y) \le L|x-y| \quad \text{for all} x,y \in \R^d \text{ for } L>0,\\
	-\infty < \underline{\tilde f} := \inf \tilde f \le \tilde f \le \sup \tilde f := \overline{\tilde f} < + \infty.
	\end{gather*}
\end{assumption}
The authors emphasize that the regularization $\tilde f$ is introduced only for technical reasons and that it does not the influence the optimization problem, as it can be shown that the dynamic stays on the hyper-surface whenever it is initialized there.

The well-posedness results for the particle and the mean-field scheme \cite{Fornasier1} read as follows.
\begin{theorem}
Let Assumption~\ref{ass:constraint} hold and $f$ with $0\le f$ be locally Lipschitz. Moreover, let $\rho_0 \in \mathcal P(\Gamma).$ For every $N\in \N$ there exists a path-wise unique strong solution $X = (X^1, \dots, X^N)$ to the system \eqref{CBOconstraint} with initial condition $X(0) = X_0.$ Moreover, it holds $X^i(t) \in \Gamma$ for all $i\in N$ and $t>0.$  
\end{theorem}

The well-posedness of the PDE is established similar to Theorem~\ref{thm:well_posed_meanfield} in Section ~\ref{sec:Variant1} with the help of an auxiliary mono-particle process $\bar X$ satisfying
\begin{align}\label{eq:mono-particle-process}
	d\bar X(t) &= -\lambda P(\bar X(t))(X(t) - v_f) dt + \sigma |\bar X(t) - v_f|P(\bar X(t)) dW(t) \nonumber\\ &\quad\qquad-\frac{\sigma^2}{2}(\bar X(t) - v_f)^2 \Delta \gamma(\bar X(t)) \nabla \gamma(\bar X(t)) dt, 
\end{align}
in strong sense for any initial data $\bar X(0) \in \Gamma$ distributed according to $\rho_0 \in \mathcal P(\Gamma).$ It holds $\text{law}(\bar X(t)) = \rho_t$ which allows to define $v_f = v_f[\rho]$ as in \eqref{eq:vfmeanfield}. For details on the well-posedness of the PDE we refer the interested reader to \cite{Fornasier1}. We proceed with the ideas leading to the rigorous mean-field limit result.

Using $N$ independent copies of this mono-particle process allows to obtain the rigorous mean-field limit with the well-known technique of Sznitman \cite{Sznitman}. In contrast to the unconstrained case, where the rigorous proof of the mean-field limit is open, the compactness of the hyper-surface makes the difference.
\begin{theorem}
Let Assumption~\ref{ass:constraint} hold and $f$ be locally Lipschitz. For any $T>0$ let $X^i(t)$ and $\bar X^i(t), i=1,\dots,N$ be solutions to $\eqref{CBOconstraint}$ or the corresponding mono-particle process, respectively, up to time $T$ with the same initial data $X^i(0) = \bar X^i(0)$ and the same Brownian motions $W^i(t)$. Then there exists a constant $C>0$ depending only on parameters, regularizations and constants, such that \[
\sup\limits_{i=1,\dots,N} \mathbb E[ |X^i(t)- \bar X^i(t)|^2] \le \frac{CT}{N}\Big(1+ CTe^{CT}\Big)
\]
holds for all $0\le t \le T.$
\end{theorem}

\begin{remark}
In addition to these results, there is a preprint \cite{Fornasier2} which reports on the convergence to the global minimizer and simulation results for applications in machine learning for the CBO scheme constrained to hyper-surfaces \eqref{CBOconstraint}. 
\end{remark}

With this we conclude the survey of the variants. In the next section we briefly discuss some applications and performance results of the variants.
 
\section{Overview of applications}\label{sec:applications}
The CBO variants were studied in various test problems. Initially, benchmark function from global optimization were used to get first results. As the variants are tailored for high dimensional applications arising in machine learning problems, they are tested against stochastic gradient descent. The preprint \cite{Fornasier2} shows some first results for the constraint method for the Ackley function on the sphere and machine learning scenarios. 

\subsection{Global Optimization problems - comparison to heuristic methods}
 In \cite{CBO1,CBOVergleich} benchmark functions from global optimization with various local minima and only one global minimum such as the Ackley, Rastrigin, Griewank, Zakharov and Wavy function were employed to test CBO against PSO and WDO. It turns out that CBO shows the best overall performance. In particular in scenarios where PSO and WDO have a very low success rate CBO leads to reasonable results with success rates $>50\%.$

\subsection{Machine Learning}
Variants 2 and 3 were tailored for applications in machine learning. A comparison between Variant 2 and the stochastic gradient descent is reported in \cite{CarrillJinZhu}. For a global optimization problem with an objective function that has many local minima the CBO variant outperforms SGD. The authors explain that this is caused by the fact that SGD needs a lot of time to escape from basins of local minima.

Another test case considers the well-known MNIST data set. Here, the differences between Variant 2 and SGD are less obvious. Nevertheless, CBO leads to slightly better results. See \cite{CarrillJinZhu} for more details.

\subsection{Global optimization with constrained state space}
The preprint \cite{Fornasier2} investigates global optimization problems from signal processing and machine learning that are naturally stated on the sphere. The first one is \textit{phase retrieval}, where the task is to recover an input vector $z$ from noisy quadratic measurements. The simulation results show that Variant 5 is able to match state-of-the-art methods for phase retrieval.

The second applications is \textit{robust subspace detection}. Here, the task is to find the principal component of a given point cloud. The performance of Variant 5 is reported to be equally good as the one of the Fast Median Subspace method applied to synthetic data. Then a computation of eigenfaces based on real-life photos from the \textit{10k US Adult Faces Database} is studied. It turns out that the results of Variant 5 are more reliable then the ones by SVD when outliers are present in the dataset.

\subsection{PDE versus SDE simulations}
In many applications of statistical physics, for example, particles in a plasma, electrostatic force or vortices in an incompressible fluid in two space dimensions, the mean-field equation is used to reduce computational cost \cite{Golse}. For consensus-based optimization we strongly recommend using the particle level for simulations. This is due to the fact, that not too many particles are needed for reasonable results and for high-dimensional problems the computation of the PDE solution is infeasible. The comparison of SDE and PDE results shown in \cite{CBO1} is just to underline  the formal limit numerically and thus to justify the analysis of the scheme on the PDE level.

\section{Conclusion, Outlook \& Open problems}
In this survey we collected the main results on Consensus-based optimization algorithms. First, we stated the original scheme on the particle level and the analytical results after a formal mean-field limit. Then we discussed variants with component-wise independent and common noise and mini-batch approaches that are tailored for high-dimensional applications arising from machine learning. A variant with component-wise common noise allows for analytical results on the particle level without passing to the limit $N\to\infty.$ 

Consensus-based optimization has similarities to the well-known Particle Swarm Optimization algorithms. Those were addressed in Section \ref{sec:PSO}, where we considered a variant that involves the personal best state of each particle in the dynamic.
The survey on the variants was completed with a section on the variants for constrained global optimization problems which involves the divergence and Laplace-Beltrami operator for hyper-planes. Then we shortly summarized some performance results of the CBO variants and mentioned comparisons to PSO, WDO and SGD. We conclude the survey with remarks on recent preprints and open problems.\\[0.5em]

\noindent\textbf{Recent preprints}\newline
This survey article discusses recent advances of the CBO model that have been published in peer-reviewed journals. Despite these, there are some preprints available which have not passed the peer-review at the time of the final version of this survey:
\begin{itemize}
\item [-] A recent preprint \cite{ADAM} proposes CBO with adaptive momentum estimation (ADAM) scheme which is well-known in the community of stochastic gradient descent methods. The article claims that the new scheme has high success rates at a low cost. Moreover, it can handle nondifferentiable activation functions in neural networks.

\item [-] As mentioned in Section~\ref{sec:PSO} there is a preprint \cite{GrassiPareschi} that discusses a SDE version of the PSO model that allows for passing to the limit $N\to\infty.$ An formal analysis on the mean-field level compares properties of CBO and PSO.

\item [-] In Section~\ref{sec:variant5} the preprint concerning the convergence to the global minimizer and machine learning application for CBO constrained to hyper-surfaces \cite{Fornasier2} was mentioned.
\end{itemize}

\noindent\textbf{Open problems}\newline
Let us mention some interesting open problems in the context of CBO: 
\begin{itemize}
	\item [-] The rigorous mean-field limit for the unconstrained method in $\R^d$ is not established. First estimates in this direction are provided in \cite{CBO2}. 
	\item [-] A convergence analysis on the particle level was only done for the component-wise common noise algorithm. A rigorous convergence analysis for other variants on the particle level remains open.
	\item [-] For comparison and qualitative performance results an estimate on the speed of convergence of the particles to the consensus-point would be of great interest. This point was mentioned in \cite{HaJinKimErrorEstimates} and is still open up to the authors knowledge.
	\item [-] In most application the structure of the objective functions is unknown, therefore one cannot guarantee the existence of a unique global minimum. This can lead to difficulties with $v_f$ for symmetric objective functions. A symmetry breaking generalization to problems with multiple global minima would therefore be very interesting.
\end{itemize}
Altogether, the analytical results and the numerical performance of the CBO variants are very promising and motivate for further research. 

\section*{acknowlegements} \noindent
CT was partly supported by the European Social Fund and by the Ministry Of Science, Research and the Arts Baden-W\"urttemberg.

\end{document}